\documentclass[a4paper]{article}
\usepackage{amsthm,amssymb,amsmath,enumerate,graphicx,epsf}

\newcommand{\COLORON}{1}
\newcommand{\NOTESON}{1}
\newcommand{\Debug}{0}
\usepackage[usenames,dvips]{color} 
\usepackage{amsthm,amssymb,amsmath,bbm,enumerate,graphicx,epsf,stmaryrd}
\usepackage[dvips, bookmarks, colorlinks=false, breaklinks=true]{hyperref} 

\newcommand{\comment}[1]{}
\newcommand{\COMMENT}[1]{}

\definecolor{darkgray}{rgb}{0.3,0.3,0.3}
\newcommand{\defi}[1]{{\color{darkgray}\emph{#1}}}

\newcommand{\acknowledgement}{\section*{Acknowledgement}}

\comment{
	\begin{lemma}\label{}	
\end{lemma}
\begin{proof}

\end{proof}

\begin{theorem}\label{}
\end{theorem} 
\begin{proof} 	

\end{proof}

}

\newtheorem{proposition}{Proposition}[section]

\newtheorem{theorem}[proposition]{Theorem}
\newtheorem{corollary}[proposition]{Corollary}

\newtheorem{lemma}[proposition]{Lemma}
\newtheorem{observation}[proposition]{Observation}
\newtheorem{conjecture}{{Conjecture}}[section]

\newtheorem{problem}[conjecture]{{Problem}}

\newtheorem{examp}[proposition]{Example}

\newcommand{\FIG}{0}

\ifnum \NOTESON = 1 \newcommand{\note}[1]{ 

\hspace*{-30pt}
	{\color{blue}  NOTE: \color{Turquoise}{\small  \tt \begin{minipage}[c]{1.1\textwidth}  #1 \end{minipage} \ignorespacesafterend }} 
	
	}
\else \newcommand{\note}[1]{} \fi

\newcommand{\afsubm}[1]{ \ifnum \Debug = 1 {\mymargin{#1}}
\fi} 

\ifnum \Debug = 1 \newcommand{\sss}{\ensuremath{\color{red} \bowtie \bowtie \bowtie\ }}
\else \newcommand{\sss}{} \fi

\ifnum \FIG = 1 
\else  \fi

\ifnum \FIG = 1 
\else  \fi

\ifnum \Debug = 1 \usepackage[notref,notcite]{showkeys}
\fi

\ifnum \COLORON = 0 \renewcommand{\color}[1]{}
\fi

\newcommand{\showFig}[2]{
   \begin{figure}[htbp]
   \centering
   \noindent
   \epsfbox{#1.eps}
   \caption{\small #2}
   \label{#1}
   \end{figure}
}

\newcommand{\N}{\ensuremath{\mathbb N}}
\newcommand{\R}{\ensuremath{\mathbb R}}
\newcommand{\Z}{\ensuremath{\mathbb Z}}

\newcommand{\cc}{\ensuremath{\mathcal C}}

\newcommand{\co}{\ensuremath{\mathcal O}}
\newcommand{\cp}{\ensuremath{\mathcal P}}

\newcommand{\cw}{\ensuremath{\mathcal W}}

\newcommand{\sig}{\ensuremath{\sigma}}

\newcommand{\one}{\mathbbm 1}
\newcommand{\sm}{\backslash}

\newcommand{\seq}[1]{\ensuremath{(#1_n)_{n\in\N}}} 

\newcommand{\g}{\ensuremath{G\ }}
\newcommand{\G}{\ensuremath{G}}

\newcommand{\Ex}{\mathbb E}
\renewcommand{\Pr}{\mathbb{P}}

\newcommand{\Tr}[1]{Theorem~\ref{#1}}

\newcommand{\Sr}[1]{Section~\ref{#1}}

\newcommand{\Prb}[1]{Problem~\ref{#1}}

\newcommand{\Or}[1]{Observation~\ref{#1}}

\newcommand{\Cg}{Cayley graph}

\newcommand{\fe}{for every}
\newcommand{\Fe}{For every}
\newcommand{\fea}{for each}

\newcommand{\st}{such that}

\newcommand{\ti}{there is}

\newcommand{\wrt}{with respect to}

\newcommand{\rw}{random walk}

\newcommand{\labtequ}[2]{
 \begin{equation} \label{#1} 	\begin{minipage}[c]{0.9\textwidth}  #2 \end{minipage} \ignorespacesafterend \end{equation} }

\newcommand{\mymargin}[1]{
  \marginpar{%
    \begin{minipage}{\marginparwidth}\small%
      \begin{flushleft}%
        {\color{blue}#1}%
      \end{flushleft}%
   \end{minipage}%
  }%
}%

\newcommand{\mySection}[2]{}

\newcommand{\rg}{random graph}
\newcommand{\gwrg}{GWRG} 
\newcommand{\PBv}{Poisson boundary}
\newcommand{\PB}{\ensuremath{\cp}}

\newcommand{\ecm}{effective conductance measure}

\title{Group-Walk Random Graphs}

\author{Agelos Georgakopoulos\thanks{Supported by EPSRC grant EP/L002787/1. This project has received funding from the European Research Council (ERC) under the European Union’s Horizon 2020 research and innovation programme (grant agreement No 639046).}\medskip 
\\
  {\small Mathematics Institute, University of Warwick, CV4 7AL, UK}}

\date{}
\begin{document}
\maketitle

\begin{abstract}
We introduce a construction that gives rise to a variety of `geometric' finite \rg  s, and describe  connections to the \PBv, Naim's kernel, and Sznitman's random interlacements.
\end{abstract}

\section{Introduction}

The purpose of this paper is to introduce a new construction of `geometric' finite \rg  s, called Group-Walk Random Graphs (\gwrg s from now on) and describe the rich connections to other objects, including the \PBv, Naim's kernel, and Sznitman's random interlacements. \gwrg s  do not only yield new interesting examples of \rg s, but as we will argue, they can be thought of as a tool for studying groups. 

\showFig{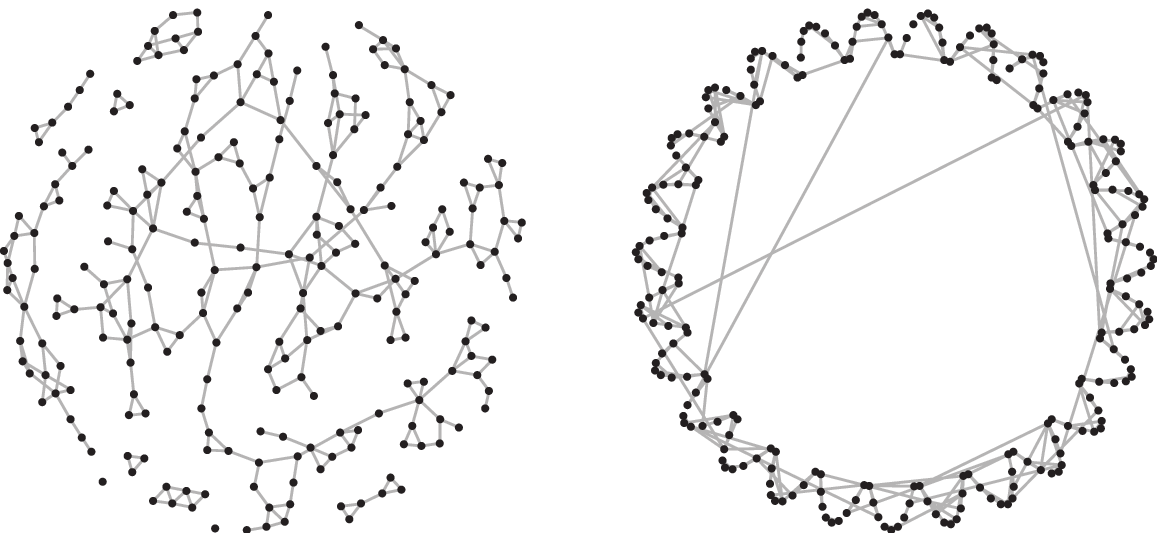}{A sample \gwrg\ produced by computer simulation by A. Janse van Rensburg. The host graph is a ternary tree, and $n=5$.}

\medskip
We start by introducing the simplest special case of \gwrg s. 
Let \g be an infinite homogeneous tree, rooted at a vertex $o$, called the \defi{host} graph; we will later allow \g to be an arbitrary locally finite \Cg, or even a more general graph.
Let $G_n:= G[ \{v\in V(G) \mid d(v,o)\leq n\}]$ be the ball of radius $n$ centered at $o$, and define the boundary $\partial G_n$ to be the set $\{v\in V(G) \mid d(v,o)= n\}$ of vertices at distance exactly $n$ from $o$. 

We construct a  \rg\ $R_n$ as follows. The vertex set of $R_n$ is the deterministic set $\partial G_n$. The edge set of  $R_n$ is constructed according to the following process. We start an independent simple \rw\ in $G_n$ from each vertex $v \in \partial G_n$, and stop it upon its first return to $\partial G_n$, letting $v^\dagger$ denote the vertex in $\partial G_n$ where this \rw\ was stopped. We then put an edge in $R_n$ joining $v$ to $v^\dagger$ \fea\ $v \in \partial G_n$. 

We could stop the construction here and declare $R_n$ to be our \gwrg, but it is more interesting to consider the following evolution: let $R^1_n:= R_n$, and for $i=2,3,\ldots$ let $R^i_n$ be the union of $R^{i-1}_n$ with an independent sample of $R^1_n$; or in other words, $R^i_n$ is the \rg\ obtained as above when we start $i$ independent particles at each vertex in $\partial G_n$.

An important observation from \cite{theta} is that these \rg s $R^i_n$ have the following scale-invariance property. Let $C,D$ be two \defi{branches} of our tree $G$, where a branch is a component of $G\sm e$ for some edge $e$. Then, for any fixed $i$, we have 

\begin{observation}[{\cite{theta}}] \label{crossings}
The expected number of edges of $R^i_n$ from branch $C$ to branch $D$ converges as $n\to \infty$. The limit is always $>0$, and it is finite if and only if $C\cap D=\emptyset$.
\end{observation}

This might at first sight look surprising, as the number of vertices of $R^i_n$ inside each of $C,D$ grows exponentially with $n$, yet no rescaling is involved in \Or{crossings}.

\showFig{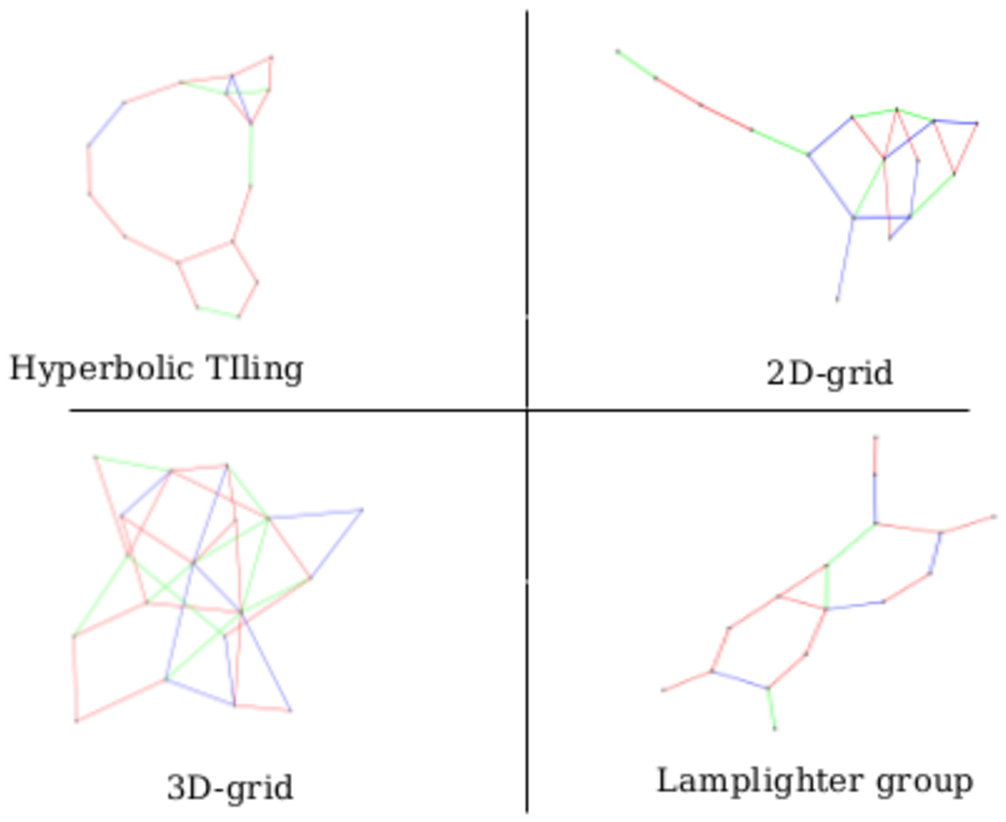}{Computer simulations of \gwrg s with 3 repetitions (coded by different colours) for different host graphs by C.~Midgley.}

The same construction can be repeated when instead of the binary tree the host graph is any infinite graph. Then Observation~\ref{crossings} has a generalisation, but in order to formulate it we need the \PBv\ \PB: instead of the `branches' of Observation~\ref{crossings} we have to talk about subgraphs `converging' to a measurable subset of \PB. This is explained in greater detail in \Sr{constr}, where we also elaborate more on the general construction of \gwrg s and their variants.
\medskip

The construction of \gwrg s was motivated by a measure space introduced in \cite{theta}, called the effective conductance measure, which is closely related to the \PBv. It is a generalisation of effective conductance for electrical networks, and it is important for the study of Dirichlet harmonic functions on infinite graphs. More details about this measure are given in \Sr{CC}.

I expect that \gwrg s can unify many existing models of geometric \rg s, while introducing new ones, and offer new tools for analysing them including the \PBv, as well as the notion of graphons \cite{LovaszLimits}. 
Conversely, \gwrg s provide an additional tool for indirectly studying groups, just like  \rw s; see \Sr{ProGr} for more.

The study of random graphs is currently one of the most active branches of graph theory. By far the most studied random graph model is that of Erd\H os \& Renyi (ER) \cite{ER}, in which every pair of vertices is joined with an edge with the same probability, and independently of each other pair. 

In recent years, many models of \defi{geometric} \rg s have been emerging \cite{HofRan,PenRgg}. The idea now is to embed the set of vertices (possibly randomly) into a geometric space ---usually the euclidean or hyperbolic plane, and their higher dimensional analogs--- and then to independently join each pair of vertices with a probability that decays as the distance between the vertices in the underlying space grows. 

One advantage of  these geometric \rg s compared to the ER model is that they can approximate real-life networks much more realistically, but they are also of great theoretical interest given the impact of the ER model. A disadvantage is that there is an infinity of such models, obtained by varying the underlying geometry, the way the points are embedded, and the connection probability as a function of distance, and no canonical choice is available. 

I like thinking of \gwrg s as geometric random graphs, where the underlying geometry is a \Cg\ $G$ of an  arbitrary finitely generated group. 
Although there is a huge variety for the underlying `geometry', the construction is in a sense canonical, and many tools for their analysis are available and are discussed here. 

\medskip
This paper is written in survey style although the material reviewed is quite new and partly under development, the main aim being to make the open problems of the project accessible to other researchers willing to get involved. A lot of the material is drawn from the paper \cite{theta} which is still in progress. New here is the definition of \gwrg s and some observations about them.

\section{The general construction of \gwrg s} \label{constr}

In the Introduction we chose the host graph $G$ to be a tree, the reason being that \Or{crossings} is easier to state in that case. Let us now
consider the general case where the host graph $G$ is arbitrary, and see how \Or{crossings} generalises, which will lead us to the definition of the \ecm.

The construction of $R^i_n$ can be repeated verbatim, except that the number of particles we start at a vertex $v\in \partial G_n$ in round $i$ is equal to the vertex degree $d_{G_n}(v)$ of $v$ in the ball $G_n$; the reason will become apparent later. However, we could instead of starting exactly $d_{G_n}(v)$ particles at $v$ in round $i$, start a random number of particles following some distribution with expectation $d_{G_n}(v)$, the most natural candidate being the Poisson distribution; the following discussion remains valid for this variant of \gwrg s. 

\Or{crossings}, which is easy to prove in the case of trees, now becomes a substantial theorem, but in order to formulate it we need to involve the concept of the Poisson boundary of $G$ to extend the above notion of branch in the correct way.

The Poisson boundary of an infinite (transient) graph \G\ is a measurable, Lebesgue-Rohlin, space $\cp=\cp(G)$, endowed with a family of probability measures $\{\mu_v \mid v\in V(G)\}$, \st\ every bounded harmonic function $h:V(G)\to \R$ can be represented by integration on $\cp$: we have $h(v)= \int_\cp \hat{h} d\mu_v(\eta)$ for a suitable boundary function $\hat{h}: \cp\to \R$.  This can be thought of as a discrete version of Poisson's integral representation formula $h(z) = \int_0^1 \frac{1-|z|^2}{|e^{2\pi i\theta} - z|^2}\hat{h}(\theta)d\theta =   \int_0^1 \hat{h}(\theta) d\nu_z(\theta)$, recovering every continuous harmonic function $h: \mathbb{D} \to \R$ in terms of its boundary values $\hat{h}: \mathbb{S}^1 \to \R$, except that we replaced $\mathbb{D}$ with a transient graph.



Triggered by the work of Furstenberg \cite{FurPoi} who introduced the concept, the study of the \PBv\ of \Cg s has grown into a very active research field; see  
\cite{erschlerICM} for a survey including many references. 
Although \ti\ a straightforward abstract construction of the \PBv\ $\cp(G)$ of any \Cg, given a concrete \g it is desirable to identify $\cp(G)$ with a geometric boundary. This pursuit can however be very hard, although some general criteria are available.

As an example, we remark that the \PBv\ of a regular tree can be identified with its set of ends, and the \PBv\ of a regular tessellation of the hyperbolic plane can be identified with its circle at infinity. More generally, the \PBv\ of any non-amenable, bounded-degree, Gromov-hyperbolic graph coincides with its hyperbolic boundary \cite{AncNeg,Ancpos}. The \PBv\ of an 1-ended bounded-degree planar graph can be identified with a circle \cite{planarPB}.

Now back to \Or{crossings}, letting \g be an arbitrary transient graph, we consider measurable subsets $X,Y$ of $\cp(G)$. One can associate with these sets sequences of vertex sets \seq{X},\seq{Y}, where $X_n,Y_n \subseteq \partial G_n$ in the above notation, \st\ for \rw\ on \g from any starting vertex, the events of converging to $X$ and visiting infinitely many of the $X_n$ coincide up to a set of measure zero, and similarly for $Y$ and the $Y_n$. The first statement of \Or{crossings} generalises: 

\begin{theorem}[\cite{theta}] \label{convC}
The expected number of edges of $R^i_n$ from $X_n$ to $Y_n$ converges as $n\to \infty$. 
\end{theorem}

For the second statement of \Or{crossings} we remark that the limit is infinite when $X\cap Y$ has positive measure, and there are, rather rare, interesting cases where the limit is infinite independently of the choice of $X,Y$ as long as they have positive measure: any lamplighter graph over a transient graph has this property \cite{theta}.

In order to understand why \Tr{convC} (or \Or{crossings}) is true, it is helpful to consider the well-known relationship between \rw s and electrical networks as introduce by Doyle \& Snell \cite{DoyleSnell}. Think of $G_n$ as an electrical network with boundary nodes $\partial G_n$, at which we impose a constant potential $v(b)=1$ for every $b\in \partial G_n$. This Dirichlet problem has the unique trivial solution $v(x)=1$ \fe\ $x\in V(G_n)$. Now, the aforementioned relationship tells us that the solution to any such Dirichlet problem can be obtained as follows: we start $d_{G_n}(b) v(b)$ \rw\ particles at each boundary node $b$, and stop them upon their first re-visit to $\partial G_n$. Then letting $v(x)$ be the expected number of visits to $x$ by all those particles divided by the degree $d(x)$ solves the Dirichlet problem \cite{script}. But as our Dirichlet problem has constant boundary values, we then expect $d(x)$ visits to each interior vertex $x$. This implies that if we start $d_{G_n}(b)$ \rw ers at each vertex $b\in \partial G_n$, stop them upon their first re-visit to $\partial G_n$, and observe the parts of their trajectories inside $G_m$ for $m<n$, then the situation we observe will be similar as if we had performed the same process on $G_m$ instead of $G_n$. This is the central observation for the proof of \Tr{convC}.

\medskip
\Tr{convC} will be crucial in the next section. 

\section{The \ecm\ \cc} \label{CC}

\Tr{convC} relates our \gwrg s to the \PBv, but in fact the connection is more intricate. Before elaborating on this, we recall
Douglas' classical formula $$E(h) = \int_0^{2\pi} \int_0^{2\pi} (\hat{h}(\eta)-\hat{h}(\zeta))^2 \Theta(\eta,\zeta) d\eta d\zeta,$$
expressing the (Dirichlet) energy of a harmonic function $h$ on the unit disc $\mathbb D$ in the complex plain from its boundary values $\hat{h}$ on the circle $\mathbb{S}^1=\partial \mathbb D$.
The physical intuition here is that $E(h)$ is the power dissipated by a circular metal plate, on the boundary circle of which some potential is imposed by an external source or field.

A discrete variant of this formula is the following, expressing the energy dissipated by a finite electrical network with a set $B$ of boundary nodes, in terms of the voltages at $B$ and certain `effective conductances' relative to $B$ (see \cite{theta} for details).
\labtequ{finEn}{$E(h) = \sum_{a,b\in B} \left(h(a)-h(b)\right)^2 C_{ab}$.}

In \cite{theta} we prove the following statement, 
providing a general formula for the Dirichlet energy of harmonic functions on a graph, which is similar to Douglas' formula

\begin{theorem}[\cite{theta}] \label{D}
\Fe\ transient graph $G$, \ti\ a measure $\cc$ on $\cp(G)^2$ \st\ \fe\ harmonic function  $h: V(G) \to \R$ with boundary function $\hat{h}$, the energy $E(h)$ equals\\ $D(\hat{h}):= \int_\cp^2 \left(\hat{h}(\eta)-\hat{h}(\zeta)\right)^2 \cc(\eta,\zeta)$. 
\end{theorem}  

This measure $\cc$, which we call the \defi{effective conductance measure}, can be thought of as a continuous analogue of effective conductance in a finite electrical network due to the similarity of the above formula with \eqref{finEn}; in the next section we will elaborate more on this.

The construction of \cc\ is based on \Tr{convC}: we set the value $\cc(X,Y)$ to be the limit value returned by that theorem, and apply Caratheodory's extension theorem to show that this defines a measure on $\PB(G)^2$. This explains the relationship between \gwrg s and \cc, and also my motivation in introducing the former.


\section{Doob's formula and Naim's kernel} \label{Doob}

\Tr{D} was motivated by a similar result of Doob \cite{Doob62}, which generalises Douglas' formula to arbitrary Green spaces. I will refrain from repeating the definition of Green spaces here, which can be found in \cite{Doob62}, and suffice it to say that they generalise Riemannian manifolds. In a sense, our \Tr{D} is a discrete version of Doob's result. But rather than working with the \PBv\ (which first appeared the year after Doob's paper \cite{FurPoi}) Doob worked with the Martin boundary, which is closely related to the \PBv\ (the latter can be defined as the support of the former, but the former is also endowed with a topology). Moreover, Doob's approach to the \ecm\ was different to ours: rather than working with the measure, Doob was working with its density, namely the Naim kernel. Let me explain this further, as it will be of interest later.

Naim  \cite{Naim} introduced the formula $$\Theta(x,y):= \frac{G(x,y)}{G(x,o) G(o,x)}$$ where $G( \cdot ,\cdot )$ denotes Green's function, and $x,y$ are points of a Green space $X$. The same formula however makes perfect sense when $x,y$ are points of a transient graph \G, in which case $G(x,y)$ denotes the expected number of visits to $y$ by random walk from $x$. 

Naim proved that $\Theta(x,y)$ can be extended to pairs of points $\{\eta,\zeta\}$ in the Martin boundary of $X$, by taking a limit 
\labtequ{limTh}{$\Theta(\eta,\zeta):= lim_{x \to \eta, y \to \zeta} \Theta(x,y).$}
The convergence of this limit is only clear when one of $x,y$ is fixed and the other converges to a point $\eta$ in the boundary, as it reduces to the well-known convergence of the Martin kernel $K(x,y):= \frac{G(x,y)}{G(o,y)}$ \cite{WoessBook09}. The two-sided convergence is a puzzling fact, and in fact is not true everywhere but `almost everywhere'\sss: the exact statement proved by Naim  \cite{Naim} is too technical to state here precisely. It involves Cartan's fine topology. As far as I know, this convergence has not been proved for graphs; we will return to this issue below.

Our intuitive interpretation of $\Theta(\eta,\zeta)$ is that it denotes the `effective conductance density' between the boundary points $\eta,\zeta$. To support this intuition, we remark in \cite{theta} that if $x,y$ are boundary vertices of a finite electrical network, then $\Theta(x,y)$ does indeed coincide with the effective conductance between $x$ and $y$ if the Green function $G$ is defined with respect to \rw\ killed at the boundary.

Doob's formula for the Dirichlet energy of a harmonic function $h$ on $X$, with boundary extension $\hat{h}$ on $M(X)$ reads
$$D(\hat{h}):= \int_{M(X)^2} \left(\hat{h}(\eta)-\hat{h}(\zeta)\right)^2 \Theta(\eta,\zeta) d\mu_o(\eta)d\mu_o(\zeta).$$
Notice the similarity to the formula of \Tr{D} and Douglas's \eqref{finEn}.

As we show in \cite{theta}, our \ecm\ \cc\ is absolutely continuous with respect to the square of harmonic measure $\mu_o$ on \PB, and so we can define a kernel $\Theta'$ on $\PB^2$ by the Radon-Nykodym derivative $\frac{\partial C}{\partial \mu_o^2} (\eta,\zeta)$. The above discussion suggests that $\Theta'$ should coincide with $\Theta$ if $\Theta$ can be defined by a limit similar to \eqref{limTh}. Now instead of trying to imitate Naim's proof of the convergence of \eqref{limTh}, which is rather technical, we propose the following

\begin{problem}
Let \g be a transient graph and $o\in V(G)$. Let \seq{x} and \seq{y} be independent simple \rw s from $o$. Then $\lim_{n,m\to \infty} \Theta(x_n,y_m)$ exists almost surely.
\end{problem}

This is seemingly easier than trying to prove the convergence of \eqref{limTh}, since for example it is easier to prove the Martingale convergence theorems than Fatou's theorem. If it is true, then one could further ask 

\begin{problem}
Let \g be a transient graph and $o\in V(G)$. Let \seq{x} and \seq{y} be independent simple \rw s from $o$, and $X,Y$ measurable subsets of $\PB(G)$. Then $\cc(X,Y)= \Ex \lim_{n,m\to \infty} \mu_{x_n}(X) \mu_{y_n}(Y) \Theta(x_n,y_m)$, where \cc\ denotes the \ecm.
\end{problem}
The factors $\mu_{x_n}(X), \mu_{y_n}(Y)$ in the above limit essentially `condition' the random walks to converge to $X,Y$ respectively.
This last problem is motivated by our expectation that $\Theta'=\Theta$.

\section{Random interlacements} \label{interl}
Recall that we defined $\cc(X,Y)$ as the limit of the expected number of edges of our \gwrg\ $R^1_n$ from $X_n$ to $Y_n$. In doing so, we were only interested in the starting and finishing vertices of our \rw s, ignoring the exact trajectories inside the host graph \G. In this section we remark that the distributions of these trajectories converge, in a certain sense, to the intensity measure of the Random Interlacement model as introduced by Sznitman \cite{SzniVac} for $\G=\Z^d$ and generalised to arbitrary transient \g by Teixera \cite{TeiInt}.

Given a  transient graph \G, the Random Interlacement on \g is  a Poisson point
process on the space of 2-way infinite trajectories in \g modulo the time shift. It is governed by a \sig-finite measure $\nu$ on $(W^*, \cw^*)$, where $W^*$ denotes the set of equivalence classes of 2-way infinite walks in \g \wrt\ the time shift, and $\cw^*$ is the canonical sigma-algebra on $W^*$. We remark that this measure $\nu$ can also be obtained from the process we used to construct our \gwrg\ $R^1_n$ as follows (the original definition of Sznitman will be given below).

Let $C$ be a cylinder set of $\cw^*$ defined by a finite walk $Z$, i.e.\ is the set of 2-way infinite walks containing $Z$ as a subsequence. Let $\cp_n$ be the random process from the definition of $R^1_n$, that is, a collection of particles performing \rw\ starting at $\partial G_n$ and stopped upon the first revisit to $\partial G_n$, with $d_{\partial G_n}(x)$ particles started at each $v\in \partial G_n$. We set 
\labtequ{mu}{$\mu(C):= \lim_n \Ex\{ \# \text{ of trajectories in $\cp_n$ containing $Z$ as a subwalk }\},$}
and extend $\mu$ to a measure on $(W^*, \cw^*)$ using e.g.\ Caratheodory's extension theorem. It turns out \cite{theta} that 
$$\mu=\nu$$

An interesting consequence of this is the following relation between  the \ecm\ $\cc$ and  the Random Interlacement intensity measure, $\nu$
\begin{corollary}
For every two measurable subsets $X,Y$ of $\PB(G)$, we have
$$\cc(X,Y)=  \nu(W^*_{XY}),$$
where  $W^*_{XY}$ denotes the set of elements of $W^*$ all initial subwalks of which meet infinitely many $X_n$ (as defined in \Sr{constr}) and all final subwalks of which meet infinitely many $Y_n$.
\end{corollary}

In order to explain why this is true we need to recall the definition of the intensity measure $\nu$ from \cite{TeiInt}. 

To begin with, given a finite subset $K$ of $V(G)$, we define the \defi{equilibrium measure} $e_K$ on $K$ by 
$$e_K(x) = \one_{x\in K} c_x \Pr_x[A \mid \text{ \rw\ does never return to $K$ }], $$
where $c_x$ is the vertex degree of $x$.

Let $ \pi^*$ be the canonical projection from the set of 2-way infinite walks $W$ to $W^*$. Then $\nu$ is defined as the unique measure on $(W^*, \cw^*)$ satisfying, for every finite subset $K$ of $V(G)$,
\labtequ{nuq}{$\one_{W^*_K} \cdot \nu = \pi^* \circ Q_K$.}
(Another way to state this formula is $\nu(\one_{W^*_K} \cdot A) = Q_K(\pi^{*-1}(A))$, where $\pi^{*-1}(A)$ returns those walks in the $\pi^*$-preimage of $A$ that enter $K$ at time 0 for the first time.)

 Here, $Q_K$ is a finite measure on the space $W_K$ of 2-way infinite walks meeting $K$, given by the formula
$$Q_K[(X_{-n})_{n\geq 0} \in A, X_0 = x, (X_{n})_{n\geq 0} \in B] = \Pr_x[A \mid \text{no return to $K$}] e_K(x) \Pr_x[B],$$
where $A,B$ are measurable subsets of the space $W_+$ of 1-way infinite walks. 
\medskip
Let me explain why $\nu$ coincides with $\mu$ as claimed above. The main idea is to think of the equilibrium measure $e_K(x)$, which is proportional to the probability of escaping $K$, as the probability for a \rw er coming from `infinity' to enter $K$ at $x$; this intuition is justified by the reversibility of our walks, and coming from infinity can be made precise using the process used to construct \gwrg s (and $\mu$). 

To make this more precise, let $G_n$ denote the ball of radius $n$ around $o$, and suppose that $K\subset V(G_n)$ for some $K\subset  V(G)$. Let $G^*_n$ be the graph obtained from $G$ by contracting the complement of $G_n$ into a single vertex $*_n$. Then the reversibility of our \rw\ implies that, for every $x \in K$, letting $\Pr_y$ denote the law of a \rw\ $X_0,X_1, \ldots X_\tau$ from $y$ stopped the first time $\tau>0$ when $X_\tau \in K \cup \{*_n\}$, we have 
$$c_x \Pr_x[X_\tau= *_n] = c_{*_n} \Pr_{*_n}[X_\tau =x] .$$
Now note that, by the  transience of \G, the limit as $n$ goes to infinity of the left hand side converges to $e_K(x)$, while the limit of the right hand side is closely related to our definition of $\mu$.


\medskip
Random interlacements have been studied extensively, and have found interesting applications in the study of the vacant set for \rw\ on discrete tori \cite{}. This connection to the \PBv\ has apparently not been observed before.

\section{Graphons}

Graphons were recently introduced \cite{LovaszLimits} as a notion of limit for sequences of dense finite graphs, and they have already had a seminal impact on combinatorics. Formally, a graphon is a symmetric, measurable function $w:[0,1]^2 \to [0,1]$. 

Every graphon naturally gives rise to a family of random graphs $G_n$ on $n$ vertices for every $n\in \N$:  we sample $n$  independent, uniformly distributed points from $[0,1]$ to be the vertices of  $G_n$, and join vertices $x,y$ with probability $w(x,y)$.

Graphs sampled this way are dense, i.e.\ have average degree of order $n$. But a variant of the above sampling method introduced in \cite{BCCZ} produces sparse graphs. 

Now note that formally, our \ecm\ \cc\ is similar to a graphon, as it is a measure on a square. In particular, we can sample a \rg\ from it. We expect these graphs to be closely related to our \gwrg s.

\begin{problem} \label{graphon}
Show that \fe\ (transient) \Cg\ \G, the sequence of corresponding \gwrg s converges to a sparse graphon in the sense of \cite{BCCZ}.
\end{problem}

\section{Simulation data} \label{sim}
Computer simulations on \gwrg\ for certain concrete host graphs $G$, performed by Chris Midgley \cite{Midgley}, suggest that certain properties are heavily influenced by  \G, while convergence with $n$ is fast enough that simulation data can help make explicit predictions. 

In most simulations the host graph was the 2 or 3-dimensional grid $\Z^2$ or $\Z^3$, the infinite binary tree $T_2$, a hyperbolic planar graph $T^\circ$ obtained from the  binary tree by adding a cycle $C_n$ joining the vertices at distance $n$ from the root \fe\ $n$, and the lamplighter graph $L$ over $\Z$.

The outcomes of these simulations, with 10,000 random graphs $R^1_n$ generated in each case for $n$ up to 8 (where the exponential growth of $T^\circ$ and $L$ become computationally demanding), suggest the following.

\begin{itemize}
 \item The number of isolated vertices of $R^1_n$ is asymptotically proportional to the number of vertices $|\partial G_n|$ of $R^1_n$. The same is true for the number of components of $R^1_n$. The corresponding leading coefficients are very similar (possibly converging to the same number) when $G$ is $\Z^2$ or $T^\circ$, slightly different when $G=L$, and very different when $G =\Z^3$.

 \item The expected diameter of the largest connected component of $R^1_n$ seems to be proportional to $\log(|\partial G_n|)$.


\end{itemize}

In a further experiment of \cite{Midgley} on the  infinite binary tree $T_2$, random graphs $R^i_n$ were generated for $2\leq n \leq 8$ for all values of $i$ until the first time that $R^i_n$ becomes connected. The data suggest (quite clearly) that the average time $i$ till connectedness is roughly $0.26 n$, i.e.\ linear in $n$. (The exact outcome was $i=(0.26 \pm 0.003)k + (0.90 \pm 0.02)$ with an adjusted $R^2$ value of 0.9993.)

\section{Further Problems}

\subsection{Dirichlet harmonic functions}

Let \g be a graph on which all harmonic functions with finite Dirichlet energy are constant; the class of such graphs is denoted by $\co_{HD}$. By \Tr{D}, we can introduce the following trichotomy for such graphs: 
\begin{enumerate}
 \item \label{ti} $\PB(G)$ is trivial (in other words, \g has the Liouville property), or
 \item \label{tii} $\PB(G)$ is not trivial, and $\cc(X,Y)=\infty$ for every two measurable $X,Y \subseteq \PB(G)$, or
\item \label{tiii} none of the above holds, but the integral of \Tr{D} is infinite for every boundary function $\hat{h}$.
\end{enumerate}

Recall that the property of being in $\co_{HD}$ is quasi-isometry invariant for graphs of bounded degree \cite{thomassenCurrents}, and therefore, given a group $\Gamma$, it is independent of the choice of the \Cg\ of $\Gamma$. A well-known open problem asks whether the Liouville property too is independent of  the choice of the  \Cg\ of $\Gamma$. The above trichotomy suggests the following refinement of that problem for groups in $\co_{HD}$

\begin{problem}
Let \g be a \Cg\ in $\co_{HD}$. Do all finitely generated \Cg s of the group of \g have the same type \wrt\ the above trichotomy?
\end{problem}

We remark that \Cg s of all three types exist: $\Z^3$ is of type \ref{ti}, lamplighter graphs over any transient \Cg\ are of type \ref{tii} \cite{theta}, while all tessellations of hyperbolic 3-space $\mathbb{H}^3$ are of type \ref{tiii}.

\subsection{\gwrg s}

None of the results suggested by the simulations of \Sr{sim} have been proved rigorously, and it would be interesting to do so, especially if the methods involved can be applied to large families of host graphs. Results on the expected number of components or isolated vertices of $R^1_n$ could be within reach.

The interesting meta-problem is to find natural properties of the \gwrg s that are universal in the sense that they are true independently of the host graph \g (under some restriction, e.g.\ \g being vertex-transitive), as well as properties that do depend on \g and the dependence can be explained. For example, what can be said about $R_n$ if the host \g is hyperbolic?

Erd\H os Renyi random graphs $G(n,p)$ are known to display a sharp  threshold for their connectedness at the value $p=\ln n/n$ for the presence of each edge. This threshold coincides with the threshold for having an isolated vertex \cite{ER}. This motivates the question of whether similar behaviour is observed for \gwrg s. To make this more precise, define the random variable $\tau_n:= \min \{ i \mid G^i_n \text{ is connected}\}$. Our first question is how concentrated $\tau_n$ is (we interpret high concentration as a sharp phase transition). The second question is whether the threshold for connectedness coincides with that for absence of isolated vertices: letting $\tau^*_n:= \min \{ i \mid G^i_n \text{ has no isolated vertices}\}$,

\begin{problem} \label{thresh}
Is $\lim_n \Ex \tau_n = \lim_n \Ex \tau^*_n$?
\end{problem}
This might depend on the host graph \G, and the answer might be positive for every `nice' host, e.g. any \Cg. 

\medskip

For any host graph \g, and $n\in \N$, consider the ball $G_n$ as a random rooted graph $G^*_n$ by rooting it at a uniformly chosen vertex in the boundary $\partial G_n$. Then, as suggested by Gourab Ray (private communication), if $G^*_n$ converges in the local weak sense, then our \gwrg\ $R^i_n$ should converge for every $i$ in the Benjamini-Schramm sense \cite{BeSchrRec}. The relationship between this limit and the host graph could be interesting. (But even if $G^*_n$ does not converge, it will have sub-sequential limits and the situation is not less interesting.)

A particular case of interest is where $G$ is the grid $\Z^d$, and we start particles at $\partial G_n$ according to the Poisson distribution in the construction of $R^i_n$. The aforementioned limit coincides then with the long-range percolation model as defined e.g.\ in \cite{NewSchulOne} (this connection was again noticed by Gourab Ray).

\comment{
Our next problem asks whether the local structure of our \gwrg s stabilises as $n$ grows
\begin{problem}
Does $R^i_n$ converge in the Benjamini-Schramm sense as $n\to \infty$ for every fixed $i$ when the host graph is vertex transitive?
\end{problem}
If the answer is positive, the limit object would provide infinite `geometric' random graphs. Moreover, the relationship between this limit and the host graph could be interesting.
}
\medskip

It would be interesting to compare our \gwrg s induced by certain simple \Cg s, like tessellations of the hyperbolic plane, to other geometric \rg s from the literature:

\begin{problem} \label{specialrg}
Show that the standard geometric \rg\ constructions can be obtained as special cases of \gwrg s by choosing an appropriate underlying graph.
\end{problem}

\subsection{Groups} \label{ProGr}

Our main objective is to understand the interplay between typical properties of \gwrg s and their host groups. In particular, we have

\begin{problem}
Which properties of the random graphs are determined by the group of the host graph and do not depend on the choice of a generating set?
\end{problem}

A further problem in this vein is

\begin{problem} \label{detect}
Are the properties of transience, Liouvilleness, existence of harmonic Dirichlet functions, and amenability on the host Cayley graph detectable by the asymptotic behaviour of the corresponding \gwrg s?
\end{problem}

Answers to these problems would make \gwrg\ a tool for studying group-theoretical questions, on a par with the study of \rw s on groups.

\section{Conclusions}
We introduced a construction of `geometric' random graphs (\gwrg s), and established strong links to the \PBv\ of their host graphs, the \ecm\ and Naim's kernel, and to random interlacements. This project will be successful if we can exploit these connection in order to make conclusions about one of these objects by studying the other. Another major aim is to relate \gwrg s to existing models of geometric random graphs. Last but not least, we would like to understand the effect of various properties of the host group on the typical graph-theoretic properties of \gwrg s.

\acknowledgement{I am grateful to Christophe Garban for suggesting the connection to Random Interlacements, which gave rise to \Sr{interl}, and to A. Janse van Rensburg and C.~Midgley for their simulations. I would like to thank Remco van der Hoffstad for triggering \Prb{thresh}.}

\bibliographystyle{plain}
\bibliography{collective}

\begin{thebibliography}{10}

\bibitem{AncNeg}
A.~Ancona.
\newblock {Negatively Curved Manifolds, Elliptic Operators, and the Martin
  Boundary}.
\newblock {\em The Annals of Mathematics}, 125(3):495, 1987.

\bibitem{Ancpos}
A.~Ancona.
\newblock Positive harmonic functions and hyperbolicity.
\newblock In {\em Potential Theory Surveys and Problems}, volume 1344 of {\em
  Lecture Notes in Mathematics}, pages 1--23. 1988.

\bibitem{BeSchrRec}
I.~Benjamini and O.~Schramm.
\newblock Recurrence of distributional limits of finite planar graphs.
\newblock {\em Electronic Journal of Probability}, 6, 2001.

\bibitem{BCCZ}
C.~Borgs, J.~T. Chayes, H.~Cohn, and Y.~Zhao.
\newblock {An {Lp} theory of sparse graph convergence I: limits, sparse random
  graph models, and power law distributions}.
\newblock Preprint 2014.

\bibitem{LovaszLimits}
C.~Borgs, J.T. Chayes, L.~Lov\'asz, V.T. S\'os, and K.~Vesztergombi.
\newblock {Convergent sequences of dense graphs I: Subgraph frequencies, metric
  properties and testing}.
\newblock {\em Advances in Mathematics}, 219(6):1801 -- 1851, 2008.

\bibitem{Doob62}
J.~L. Doob.
\newblock {Boundary properties of functions with finite Dirichlet integrals.}
\newblock {\em Ann. Inst. Fourier}, 12:573--621, 1962.

\bibitem{DoyleSnell}
P.~G. Doyle and J.~L. Snell.
\newblock {\em Random Walks and Electrical Networks}.
\newblock Carus Mathematical Monographs 22, Mathematical Association of
  America, 1984.

\bibitem{ER}
P.~Erd\"os and A.~R\'enyi.
\newblock {On random graphs I.}
\newblock {\em Publ. Math. Debrecen}, 6:290--297, 1959.

\bibitem{erschlerICM}
A.~Erschler.
\newblock {Poisson-Furstenberg Boundaries, Large-scale Geometry and Growth of
  Groups}.
\newblock In {\em Proceedings of the ICM}, pages 681--704, 2010.

\bibitem{FurPoi}
H.~Furstenberg.
\newblock {{A Poisson formula for semi-simple Lie groups.}}
\newblock {\em The Annals of Mathematics}, 77(2):335--386, 1963.

\bibitem{planarPB}
A.~Georgakopoulos.
\newblock The boundary of a square tiling of a graph coincides with the poisson
  boundary.
\newblock To appear in {\em Invent.\ Math.}, DOI 10.1007/s00222-015-0601-0.

\bibitem{script}
A.~Georgakopoulos.
\newblock Electrical networks from the mathematical viewpoint.
\newblock Lecture notes. {I}n preparation.

\bibitem{theta}
A.~Georgakopoulos and V.~Kaimanovich.
\newblock In preparation.

\bibitem{HofRan}
Remco Van~Der Hofstad.
\newblock Random graphs and complex networks.
\newblock Lecture Notes, 2013.

\bibitem{NewSchulOne}
Newman~C. M. and Schulman~L. S.
\newblock {One Dimensional $1/|j − i|s$ Percolation Models: The Existence of
  a Transition for $s \leq 2$}.
\newblock {\em Commun.\ Math.\ Phys.}, 104:547--571, 1986.

\bibitem{Midgley}
C.~Midgley.
\newblock Random graphs from groups, 2014.
\newblock {Undergraduate research project. University of Warwick}.

\bibitem{Naim}
L.~Na\"im.
\newblock Sur le r\^ole de la fronti\`ere de {R.\ S.\ Martin} dans la th\'eorie
  du potentiel.
\newblock {\em Annales Inst.\ Fourier}, 7:183--281, 1957.

\bibitem{PenRgg}
Mathew Penrose.
\newblock {\em Random Geometric Graphs}.
\newblock Oxford University Press, 2003.

\bibitem{SzniVac}
A.-S. Sznitman.
\newblock Vacant set of random interlacements and percolation.
\newblock {\em Annals of Mathematics}, 171(3):2039--2087, 2010.

\bibitem{TeiInt}
Augusto Teixeira.
\newblock Interlacement percolation on transient weighted graphs.
\newblock 14:1604--1627, 2009.

\bibitem{thomassenCurrents}
C.~Thomassen.
\newblock Resistances and currents in infinite electrical networks.
\newblock {\em J.~Combin.\ Theory (Series B)}, 49:87--102, 1990.

\bibitem{WoessBook09}
Wolfgang Woess.
\newblock {\em {Denumerable Markov chains. Generating functions, boundary
  theory, random walks on trees.}}
\newblock {EMS Textbooks in Mathematics. European Mathematical Society (EMS),
  Z\"urich}, 2009.

\end{thebibliography}
\end{document}